\documentstyle{article}

\let\s=\sigma
\let\o=\omega

\def\E{\mathop{\hbox{\sf E}}\nolimits}
\def\P{\mathop{\hbox{\sf P}}\nolimits}

\def\RR{{\cal R}}
\def\KK{{\cal K}}

\let\o=\omega

\let\s=\sigma
\let\o=\omega

\let\detla=\delta

\title{Local Energy Statistics in Directed Polymers.}
\author{Irina Kurkova
\thanks{Laboratoire de Probabilit\'es et Mod\`eles Al\'eatoires,
  Universit\'e Paris 6, B.C. 188; 4, place Jussieu, 75252 Paris Cedex
  05, France. E-mail : kourkova@ccr.jussieu.fr }}
\date{}

\begin{document}

\maketitle

\noindent{\bf Abstract}
  Recently, Bauke and Mertens conjectured that the local statistics
    of energies in random spin systems with discrete
      spin space should, in most circumstances,
        be the same as in the random energy model.
     We show that this conjecture holds true as well
    for directed polymers in random environment.
   We also show that, under certain conditions,
     this conjecture holds for directed polymers
       even if energy levels that grow moderately
        with the volume of the system are considered.

\bigskip


\newtheorem{theo}{Theorem}
\newtheorem{cor}{Corollaire}
\newtheorem{conj}{Conjecture}
\newtheorem{lem}{Lemme}
\newtheorem{defi}{Definition}
\newtheorem{exe}{Exercice}
\newtheorem{pro}{Proposition}

\section{Introduction and results}
      Recently, Bauke and Mertens have proposed in \cite{BaMe} a
        new and original  look at disordered spin systems.
This point of view consists of studying the micro-canonical
scenario,
  contrary to the canonical formalism, that has become
  the favorite tool to treat  models of statistical mechanics.
  More precisely, they analyze the statistics of spin
  configurations whose energy is very close to a given value.
  In discrete spin systems, for a given system size,
  the Hamiltonian will take on a finite number
  of random values, and generally
  (at least, if the disorder is continuous)
        a given value $E$ is attained
     with probability $0$.
   One may, however, ask :
   How close to $E$ the best approximant is
      when the system size grows and, more generally,
  what the distribution of the energies that come closest to
  $E$ is~? Finally, how the values of the corresponding
  configurations are distributed in configuration space~?

   The original motivation for this viewpoint came from
   a reformulation of a problem in combinatorial optimization,
   the number partitioning problem
   (this is the problem of
    partitioning $N$ (random) numbers into two subsets such that
    their sums in these subsets are as close as possible)
          in terms of a spin system
    Hamiltonian \cite{BFM, M1,M2}. Mertens conjecture stated in these
      papers has been proven  to be correct in \cite{BCP}
        (see also \cite{BCMP}),
    and generalized in \cite{BK1}
        for the partitioning into $k>2$ subsets.

     Some time later, Bauke and Mertens generalized this conjecture
       in the following sense : let $(H_N(\sigma))_{\s \in \Sigma_N}$
         be the Hamiltonian
   of any disordered spin system with discrete spins
      ($\Sigma_N$ being the configuration space)  and
   continuously  distributed couplings, let $E$ be any given number,
   then the distribution of the close to optimal approximants of the
     level $\sqrt{N}E$ is asymptotically
       (when the volume of the system $N$ grows to infinity)
             the same as if the energies $H_N(\s)$
       are replaced by independent Gaussian random
        variables with the same mean and variance as $H_N(\s)$
     (that is the same as for Derrida's Random Energy spin glass Model \cite{D1},
        that is why it is called the REM conjecture).

    What this distribution for independent Gaussian random variables
    is ? Let $X$ be a standard Gaussian random variable,
      let $\delta_N \to 0$ as $N \to \infty$, $E \in {\bf R}$,
          $b>0$.
       Then it is easy to compute that
       $$ \P(   X  \in [E -\detla_N b, E+ \delta_N b])= (2 \delta_N b)
       \sqrt{1/(2\pi)}e^{-E^2/2}(1+o(1))\ \ \ N \to \infty.$$
   Let now $(X_\s)_{s \in \Sigma_N}$ be  $|\Sigma_N|$
   independent standard Gaussian random variables.
    Since they are independent, the number of them
    that are in the interval $[E -\detla_N b, E+ \delta_N
     b]$ has a Binomial distribution with parameters
   $(2 \delta_N b)
       \sqrt{1/(2\pi)}e^{-E^2/2}(1+o(1))$ and $ |\Sigma_N|$. If we put
          $$\delta_N =|\Sigma_N|^{-1} \sqrt{2\pi}(1/2)e^{E^2/2},$$
   by a well known theorem of the course of elementary Probability,
   this random number converges in law  to the Poisson
     distribution with  parameter $b$ as $N \to \infty$. More generally,
  the point process
   $$ \sum_{\sigma \in \Sigma_N} \delta_{ \{\delta_N^{-1} N^{-1/2}
        |\sqrt{N} X_\s-
     \sqrt{N}E|\}} $$
 converges, as $N \to \infty$, to the Poisson point process
      in ${\bf R}_+$ whose intensity measure is the Lebesgue measure.

  So,  Bauke and Mertens conjecture states that
    for the Hamiltonian $(H_N(\s))_{\s \in \Sigma_N}$
    of any disordered spin system and
       for a suitable normalization $C(N,E)$
  the sequence of point processes
  $$ \sum_{\sigma \in \Sigma_N} \delta_{ \{C(N,E)|H_N(\s)-
     \sqrt{N}E|\}} $$
      converges, as $N \to \infty$, to the Poisson point process
      in ${\bf R}_+$ whose intensity measure is the Lebesgue measure.
           In other words, the best approximant
     to $\sqrt{N} E$ is at distance $C^{-1}(N,E)W$,
      where $W$ is an exponential random variable of mean $1$.
    More generally, the  $k$th best approximant
     to $\sqrt{N} E$ is at distance $C^{-1}(N,E)(W_1+\cdots +W_k)$,
       where $W_1,\ldots, W_k$ are  independent
       exponential random variables of mean $1$, $k=1,2\ldots$
            It appears rather surprising that such a result
              holds in great generality.
      Indeed, it is well known that the correlations of the random
        variables are strong enough to modify e.g.\ the maxima of the
        Hamiltonian.
    This conjecture
      has been proven in \cite{BK2}
   for a rather large class of disordered spin systems
     including
     short range lattice spin systems as well as
    mean-field  spin glasses,
        like $p$-spin Sherringthon-Kirkpatrick (SK) models with
          Hamiltonian $H_N(\s)=N^{1/2-p/2} \sum_{i_1,\ldots, i_p}
             \s_{i_1}\cdots \s_{i_p}J_{1\leq i_1,\ldots, i_p\leq N}$
               where $J_{i_1,\ldots, i_p}$ are
        independent standard Gaussian random variables, $p\geq 1$.
              See also \cite{BCMN1}
     for the detailed study of the case $p=1$.

   Two questions naturally pose themselves.
    (i) Consider instead of $E$, $N$-dependent
      energy levels, say, $E_N={\rm const} N^\alpha$.
        How fast can we allow $E_N$ to grow with $N \to \infty$
          for the same behaviour
 (i.e.\ convergence to the standard Poisson point process under  a
  suitable  normalization) to hold ?
     (ii) What type of behaviour can we expect
        once $E_N$ grows faster than this value ?

    The first question (i) has been investigated
    for Gaussian disordered spin systems in \cite{BK2}.
           It turned out that for short range
     lattice spin systems on ${\bf Z}^d$  this
     convergence is still true up to $\alpha<1/4$.
              For mean-field  spin glasses,
        like $p$-spin SK models with
          Hamiltonian $H_N(\s)=N^{1/2-p/2} \sum_{i_1,\ldots, i_p}
             \s_{i_1}\cdots \s_{i_p}J_{i_1,\ldots, i_p}$
               mentioned above,
     this conjecture holds true up to $\alpha<1/4$
       for $p=1$ and up to $\alpha<1/2$ for $p\geq 2$.
             It has been proven in \cite{BCMN2}
               that the conjecture fails at
                 $\alpha=1/4$ for $p=1$ and $\alpha=1/2$
                   for $p=2$.
              The paper \cite{BCMN2}
              extends also these results for non-Gaussian
   mean-field $1$-spin SK models with $\alpha>0$.

  The second question (ii), that is the local behaviour
      beyond  the critical value of $\alpha$,
          where Bauke and Mertens conjecture fails,
               has been investigated
      for Derrida's Generalized Random
   Energy Models  (\cite{D2}) in \cite{BK3}.

    Finally,
    the paper \cite{BGK} introduces a new REM conjecture,
    where the range of energies involved is not reduced to a small
    window. The authors prove that for  large class of random Hamiltonians
     the point process of properly normalized energies
       restricted to a sparse enough random subset of spin
     configuration space converges to the same point process
       as for the Random Energy Model, i.e. Poisson point process
      with intensity measure $\pi^{-1/2}e^{-t\sqrt{2\ln 2}}dt$.

      In this paper we study Bauke and Merten's conjecture
          on the local behaviour of energies not
      for  disordered spin systems but for directed
      polymers in random environment.
       These models have received enough of attention
       of mathematical community over past fifteen years,
        see e.g.\ \cite{CSY} for a survey of the main results
        and references therein.
           Let $(\{w_n\}_{n\geq 0}, P)$ is  a simple
             random walk on the $d$-dimensional lattice
               ${\bf Z}^d$. More precisely,
               we let $\Omega$ be the path space
$\Omega=\{\omega=(\omega_n)_{n\geq 0};
  \omega_n\in {\bf Z}^d, n\geq 0\}$,
    ${\cal F}$ be the cylindrical $\sigma$-field on $\Omega$
     and for all $n\geq 0$, $\omega_n: \omega \to \omega_n$
       be the projection map.
         We consider the unique probability measure $P$
         on $(\Omega, {\cal F})$ such that $\omega_1-\omega_0,
         \ldots, \omega_n-\o_{n-1}$ are independent and
         $$ P(\o_0=0)=1,\ \
           P(\o_n-\o_{n-1}=\pm \delta_j)=(2d)^{-1}, \ \
             j=1,\ldots, d,$$
 where $\delta_j=(\delta_{kj})_{k=1}^d$ is the $j$th
    vector of the canonical basis of ${\bf Z}^d$.
       We will denote by $S_N=\{\omega^N=(i,\omega_i)_{i=0}^N\}$
          ($(i,\omega_i)\in {\bf N}\times {\bf Z}^d$)
         the space of paths of length $N$.
   We define the energy of the path $\omega^N=(i,\omega_i)_{i=0}^N$
   as
\begin{equation}
\label{enn}
     \eta(\omega^N)=N^{-1/2}\sum_{i=1}^N \eta(i,\o_i)
     \end{equation}
where $\{\eta(n,x) : n \in {\bf N}, x\in {\bf Z}^d\}$
   is a sequence of independent identically distributed
     random variables  on a probability space $(H, {\cal G}, \P)$.
  We assume that they have  mean zero  and variance $1$.

Our first theorem  extends Bauke and Merens conjecture
  for directed polymers.

\begin{theo}
\label{th0}
 Let $\eta(n,x)$, $\{\eta(n,x) : n \in {\bf N}, x\in {\bf
Z}^d\}$,
   be the i.i.d. random variables of the  third moment finite and
     with the Fourier transform
     $\phi(t)$ such that $|\phi(t)|=O(|t|^{-1})$, $|t|\to \infty$.
        Let $E_N=c \in {\bf R}$ and let
 \begin{equation}
     \label{delta0}
      \delta_N = \sqrt{\pi/2} e^{c^2/2}
           ((2d)^N)^{-1}.
      \end{equation}
    Then the point process
\begin{equation}
\label{th1e0}
  \sum_{\o^N\in S_N} \delta_{\{\delta_N^{-1}
  |\eta(\o^N)-E_N|\}}
\end{equation}
   converges weakly as $N \uparrow \infty$ to the Poisson
      point process ${\cal P}$ on ${\bf R}_+$
        whose intensity measure is the Lebesgue measure.
    Moreover, for any $\epsilon>0$ and any $b \in {\bf R}_+$
\begin{equation}
\label{th1b0} \P(\forall N_0\ \exists N \geq N_0,\
 \exists \o^{N,1}, \o^{N,2}\ : \
   {\rm cov}\,(\eta(\o^{N,1}), \eta(\o^{N,2}))>\epsilon\ :$$
    $$  |\eta(\o^{N,1})-E_N|\leq  |\eta(\o^{N,2})-E_N|\leq \delta_N
      b)=0.
\end{equation}
\end{theo}

  The decay assumption on the Fourier transform is not optimal,
  we believe that it can be weaken but we did not try to
    optimize it. Nevertheless, some condition
  of this type is needed, the result can not be extended
  for discrete distributions where the number of possible
  values the Hamiltonian takes on would be finite.

The next two theorems
     prove Bauke and Mertens conjecture
  for directed polymers in Gaussian environment for growing levels
  $E_N=cN^{\alpha}$.
    We are able to prove that this conjecture
      holds true for $\alpha<1/4$ for polymers
        in dimension $d=1$ et
      and $\alpha<1/2$ in dimension
        $d\geq 2$.
          We leave this investigation
         open for non-Gaussian environments.

     The values $\alpha=1/4$ for $d=1$ and
     $\alpha=1/2$ for $d\geq 2$ are likely to be the true
     critical values. Note that these are the same
    as for Gaussian SK-spin glass models
    for $p=1$ and $p=2$ respectively according to
     \cite{BCMN2}, and likely for $p\geq 3$ as well.

\begin{theo}
\label{th1}
 Let $\eta(n,x)$, $\{\eta(n,x) : n \in {\bf N}, x\in {\bf
Z}^d\}$, be independent standard Gaussian random variables.
 Let $d=1$. Let $E_N=c N^{\alpha}$ with
  $c \in {\bf R}$, $\alpha \in [0, 1/4[$ and
     \begin{equation}
     \label{delta}
      \delta_N = \sqrt{\pi/2} e^{E_N^2/2}
           (2^N)^{-1}.
           \end{equation}
   Then the point process
\begin{equation}
\label{th1e}
  \sum_{\o^N\in S_N} \delta_{\{\delta_N^{-1}
  |\eta(\o^N)-E_N|\}}
\end{equation}
   converges weakly as $N \uparrow \infty$ to the Poisson
      point process ${\cal P}$ on ${\bf R}_+$
        whose intensity measure is the Lebesgue measure.
    Moreover, for any $\epsilon>0$ and any $b \in {\bf R}_+$
\begin{equation}
\label{th1b} \P(\forall N_0\ \exists N \geq N_0,\
 \exists \o^{N,1}, \o^{N,2}\ : \
   {\rm cov}\,(\eta(\o^{N,1}), \eta(\o^{N,2}))>\epsilon\ :$$
    $$  |\eta(\o^{N,1})-E_N|\leq  |\eta(\o^{N,2})-E_N|\leq \delta_N
      b)=0.
\end{equation}
\end{theo}

\begin{theo}
\label{th2}
Let $\eta(n,x)$, $\{\eta(n,x) : n \in {\bf N}, x\in {\bf
Z}^d\}$ be independent standard Gaussian random variables.
 Let $d \geq 2$. Let $E_N=c N^{\alpha}$ with
  $c \in {\bf R}$, $\alpha \in [0, 1/2[$  and
\begin{equation}
\label{delta1}
      \delta_N = \sqrt{\pi/2} e^{E_N^2/2}
           ((2d)^N)^{-1}.
           \end{equation}
   Then the point process
\begin{equation}
\label{th2e}
  \sum_{\o^N\in S_N} \delta_{\{\delta_N^{-1}
  |\eta(\o^N)-E_N|\}}
\end{equation}
   converges weakly as $N \uparrow \infty$ to the Poisson
      point process ${\cal P}$ on ${\bf R}_+$
        whose intensity measure is the Lebesgue measure.
 Moreover, for any $\epsilon>0$ and any $b \in {\bf R}_+$
\begin{equation}
\label{th2b} \P(\forall N_0\ \exists N \geq N_0,\
 \exists \o^{N,1}, \o^{N,2}\ : \
   {\rm cov}\,(\eta(\o^{N,1}), \eta(\o^{N,2}))>\epsilon\ :$$
    $$  |\eta(\o^{N,1})-E_N|\leq  |\eta(\o^{N,2})-E_N|\leq \delta_N
      b)=0.
\end{equation}
\end{theo}

\noindent{\bf Acknowledgements.}
  The author thanks Francis Comets for introducing him
 to the area of directed polymers. He also thanks
  Stephan Mertens and Anton Bovier for attracting
  his attention to the local behavior of disordered spin systems
 and interesting discussions.

\section{Proofs of the theorems.}

Our approach  is based on the following  sufficient condition
   of convergence to the Poisson point process.
   It has been proven in a somewhat more general form
   in \cite{BK1}.

\begin{theo}
\label{tc}
 Let $V_{i,M}\geq 0$, $i\in {\bf N}$, be a family of
non-negative random
  variables satisfying the following assumptions : for any
   $l \in {\bf N}$ and all sets of constants $b_j>0$,
   $j=1,\ldots,l$
   $$ \lim_{ M \to \infty} \sum_{(i_1,\ldots, i_l) \in \{1,\ldots,
   M\} }\P(\forall_{j=1}^{l} V_{i_j, M}<b_j) \to \prod_{j=1}^l b_j $$
   where the sum is taken over all possible sequences of
   \underline{different}
    indices $(i_1,\ldots, i_l)$.
  Then the point process
  $$ \sum_{i=1}^M \delta_{\{V_{i,M}\}}$$
  on ${\bf R}_+$ converges weakly in distribution as
  $M \to \infty$
to the Poisson
      point process ${\cal P}$ on ${\bf R}_+$
        whose intensity measure is the Lebesgue measure.
\end{theo}

  Hence, in  all our proofs,
   we just have to verify the hypothesis of Theorem \ref{tc}
        for $V_{i,M}$ given by
  $\delta_N^{-1}
  |\eta(\o^{N,i})-E_N|$, i.e.\ we must show that
\begin{equation}
\label{zet}
   \sum_{(\o^{N,1}, \ldots, \o^{N,l}) \in S_N^{\otimes l}}
       \P(\forall_{i=1}^{l} :  |\eta(\o^{N,i})-E_N|<b_i \delta_N)
             \to b_1\cdots b_l
 \end{equation}
   where the sum is taken over all sets of different paths
$(\o^{N,1}, \ldots, \o^{N,l})$.

\bigskip

\noindent{\bf Informal proof of Theorem \ref{th0}.}
  Before proceeding with rigorous proofs let us give some
  informal arguments supporting  Theorem~\ref{th0}.

The random variables $\eta(\o^{N,i})$, $i=1,\ldots,l$,
  are the sums of independent identically distributed
  random variables with mean $0$ and the covariance matrix
   $B_N(\o^{N,1},\ldots, \o^{N,l})$ with $1$ on the diagonal
    and the covariances
       ${\rm cov}\, (\eta(\o^{N,i}),
         \eta(\o^{N,j}))= N^{-1}\#\{m : \o^{N,i}_m=\o^{N,j}_m\}
          \equiv b_{i,j}(N)$.

   The number of sets $(\o^{N,1},\ldots, \o^{N,l})$
   with $b_{i,j}(N)=o(1)$ ($o(1)$ should be chosen of an appropriate order)
        for all pairs $i\ne j$,
    $i,j=1,\ldots,l$, as $N\to \infty$, is $(2d)^{Nl}(1-\gamma(N))$
      as $N \to \infty$ where $\gamma(N)$ is exponentially small
        in $N$.
      For all such sets
        $(\o^{N,1},\ldots, \o^{N,i})$,
          by the local Central Limit Theorem,
    the random variables $\eta(\o^{N,i})$, $i=1,\ldots,l$,
    should behave asymptotically as
    Gaussian random variables with covariances
      $b_{i,j}(N)=o(1)$ and the determinant of the covariance matrix
   $1+o(1)$. Therefore, the probability
         that these random variables belong to $[-\delta_N b_i+c, \delta_N
    b_i+c]$ respectively for $i=1,\ldots,l$,
       equals $$(2\delta_N b_1)\cdots (2\delta_N b_l)
       (\sqrt{2\pi})^{-l} e^{-c^2l/2}=
          b_1\cdots b_l 2^{-Nl}(1+o(1)).$$
          Since the number of such sets
         $(\o^{N,1},\ldots, \o^{N,l})$
            is $(2d)^{Nl}(1+o(1))$,
   the sum (\ref{zet}) over them converges to $b_1\cdots b_l$.

   Let us turn to the remaining tiny part of $S_N^{\otimes l}$
     where $(\o^{N,1},\ldots, \o^{N,l})$
  are such that the covariances
     $b_{i,j}(N)\ne o(1)$ with $o(1)$ of an appropriate order
         for some $i\ne j$,
   $i,j=1,\ldots, l$, $N\to \infty$.
   The number of such sets is exponentially smaller than
     $(2d)^{Nl}$.
       Here two possibilities should be
considered differently.

 The first one is when the covariance matrix
  is non-degenerate. Then, invoking again the Central Limit Theorem,
    the probabilities $\P(\cdot)$
      in this case are not greater than
      $$ ({\rm det}B_N(\o^{N,1},\ldots, \o^{N,l}))^{-1/2}
           (2\delta_N b_1)\cdots(2\delta_N b_l)
       (\sqrt{2\pi})^{-l}.$$
    From the definition of the covariances
    of $\eta(\o^{N,i})$,
    ${\rm det}\,B_N(\o^{N,1},\ldots, \o^{N,l})$ is a finite polynomial in
    the variables $1/N$. Therefore the probabilities $\P(\cdot)$ are bounded
      by $(2d)^{-Nl}$ up to a polynomial term, while
   the number of sets $(\o^{N,1},\ldots, \o^{N,l})$
    such that $b_{i,j}(N)\ne o(1)$ some $i\ne j$,
   $i,j=1,\ldots, l$, is exponentially smaller than
     $(2d)^{Nl}$. Therefore the sum (\ref{zet})
   over such sets $(\o^{N,1},\ldots, \o^{N,l})$
     converges to zero exponentially fast.

        Let now $(\o^{N,1},\ldots, \o^{N,l})$
          be such that $B_N(\o^{N,1},\ldots, \o^{N,l})$
   is degenerate of the rank $r<l$.
     Then, without loss of generality, we may assume
         that $\eta(\o^{N,1}),\ldots, \eta(\o^{N,r})$
     are linearly independent, while
   $\eta(\o^{N,r+1}),\ldots, \eta(\o^{N,l})$
   are their linear combinations.
        Then the  probabilities $\P(\cdot)$
  are bounded by the probabilities
  that only $\eta(\o^{N,1}),\ldots, \eta(\o^{N,r})$
   belong to the corresponding  intervals,
  which are at most $2^{-Nr}$ up to a polynomial term as
   previously.
      Moreover, we will show that for no one $m=0,1,\ldots, N$,
  $\o^{N,1}_m,\ldots, \o^{N,r}_m$ can not be all different.
     Otherwise, each of $\o^{N,r+1},\ldots, \o^{N,l}$
       would coincide with one of
    $\o^{N,1},\ldots, \o^{N,r}$, which is impossible since the
      sum (\ref{zet}) is taken over sets of different(!) paths.
     This implies that the number of such sets
       $(\o^{N,1},\ldots, \o^{N,r})$
   is exponentially smaller than $2^{Nr}$. Furthermore,
     the number of possibilities to complete each of these sets
       by  $\o^{N,r+1},\ldots, \o^{N,l}$
   such that  $\eta(\o^{N,r+1}),\ldots, \eta(\o^{N,l})$
   are linear combinations
 of     $\eta(\o^{N,1}),\ldots, \eta(\o^{N,r})$
         is $N$-independent. Thus the number of sets
   $(\o^{N,1},\ldots, \o^{N,l})$ in this case being exponentially
   smaller than $2^{Nr}$, and the probabilities
      being $2^{-Nr}$ up to a polynomial term,
          the corresponding sum (\ref{zet})
   converges to zero. This completes the informal proof
   of (\ref{th1e0}) in Theorem \ref{th0}.

 We now give rigorous proofs. We start with proofs of Theorems
  \ref{th1} and \ref{th2} in Gaussian environment
  and give the proof of Theorem \ref{th0} after that.

\bigskip

\noindent{\bf Proof of Theorem~\ref{th1}.}
  For  $\eta \in ]0, 1/2[$  let us denote by
\begin{equation}
\label{rr}
  {\cal R}_{N,l}^{\eta}=\{ (\o^{N,1},\ldots, \o^{N,l}) :
    \ {\rm cov}(\eta(\o^{N,i}),\eta( \o^{N,j}))\leq N^{\eta-1/2},\
    \forall i,j=1,\ldots,l,\
      i\ne j\}.
      \end{equation}

 \noindent{\it Step 1.} As a first preparatory step, we need to estimate the capacity
     of ${\cal R}_{N,l}^{\eta}$ in (\ref{zgu}).
    Let us first note that for any two paths
     $\o^{N,1}, \o^{N,2} \in S_N$
$${\rm cov}(\eta(\o^{N,1}), \eta(\o^{N,2}))=s/N$$ if and only if
  $$ \#\{m : (\omega_m^1, m)=(\omega_m^2,m) \}=s,$$
    i.e. the number of moments of time within the period
    $[0,N]$ when the trajectories $\o^{N,1}$ and $\o^{N,2}$
       are at the same point of the space ${\bf Z}$ equals $s$.
      But due to the symmetry of the simple random walk
     \begin{equation}
       \label{rr}
       \#\Big\{\omega^{N,1},\omega^{N,2} :
                  \#\{m \in[0,\ldots,N] :
                \omega_m^1 - \o_m^2=0\} =s \Big\}$$
       $$ = \#\Big\{\omega^{N,1},\omega^{N,2} :
                   \#\{m \in[0,\ldots,N] :
             \omega_m^1 + \o_m^2=0\} =s \Big\}.
             \end{equation}
  Taking into account the fact that the random walk starting from $0$
    can not visit $0$ at odd moments of time, we obtain that
      (\ref{rr})
      equals
        $$ \#\Big\{\omega^{2N} :
                    \#\{m \in[0,\ldots,2N] :
                \omega_m = 0\}=s \Big\}.$$
  This last number is well-known for the simple random walk
      on ${\bf Z}$ : it equals
     $ 2^{2N} 2^{s-2(2N)}{2N \choose 2(2N)-s}$
       (see e.g. \cite{Fe})
         which is, by Stirling's formula,
                  when $s=[N^{1/2+\eta}]$,
          $\eta \in ]0, 1/2[$, equivalent to
     $2^{2N} (2\pi N)^{-1/2} e^{-s^2/(2(2N))}=
        2^{2N} (2\pi N)^{-1/2} e^{-N^{2\eta}/4}$ as $N \to\infty$.
        Finally, we obtain that for all $N \geq 0$
            the number  (\ref{rr})
          it is not greater than $2^{2N}e^{-h N^{2\eta}}$
            with some constant $h>0$.
 It follows that for all $N>0$
\begin{eqnarray}
 \lefteqn{ |S_N^{\otimes,l}\setminus {\cal R}_{N,l}^{\eta}|
 }\nonumber \\
      &\leq & (l(l-1)/2) 2^{N(l-2)}
        \#\Big\{\omega^{N,1},\omega^{N,2} :
                  \#\{m \in[0,\ldots,N] :
                \omega_m^1 - \o_m^2=0\} \geq N^{1/2+\eta}\Big\}
                \nonumber \\
    &\leq & 2^{Nl} C N \exp(-h N^{2\eta}) \label{zgu}
\end{eqnarray}
   where $C>0$, $h>0$ are some constants.

\medskip

\noindent{\it Step 2.} The second preparatory step
     is the estimation (\ref{es1}) and (\ref{es2})
        of the probabilities in the sum (\ref{zet}).
 Let $B_N(\o^{N,1},\ldots, \o^{N,l})$
    be the covariance matrix of the random variables
        $\eta(\o^{N,i})$  for
    $i=1,\ldots, l$.
    Then, if $B_N(\o^{N,1},\ldots, \o^{N,l})$ is non-degenerate,
    \begin{equation}
    \label{mia}
 \P(\forall_{i=1}^{l} : |\eta(\o^{N,i})-E_N|<b_i \delta_N)
 = \int_{C(E_N)} \frac{e^{-(\vec z B_N^{-1}(\o^{N,1},\ldots,
 \o^{N,l})\vec z)/2}}{(2\pi)^{l/2}
   \sqrt{{\rm det} B_N(\o^{N,1},\ldots,
 \o^{N,l})} }\, d\vec z
 \end{equation}
   where
   $$ C(E_N)=\{\vec z=(z_1,\ldots, z_l) : |z_i- E_N|\leq
   \delta_N b_i, \forall i=1,\ldots,l \}.$$

     Let $\eta \in ]0, 1/2[$.
      Since $\delta_N$ is exponentially small in $N$, we see that
      uniformly for
        $(\omega^{N,1},\ldots, \o^{N,l}) \in {\cal
       R}_{N,l}^{\eta}$,
  the probability (\ref{mia}) equals
   \begin{equation}
   \label{es1}
   (2\delta_N/\sqrt{2\pi})^l (b_1\cdots b_l)
e^{-(\vec E_N B_N^{-1}(\o^{N,1},\ldots,
 \o^{N,l})\vec E_N)/2}(1+o(1))$$
$$ =(2\delta_N/\sqrt{2\pi})^l (b_1\cdots b_l)
   e^{-\|E_N\|^2(1+O(N^{\eta-1/2}))/2}(1+o(1))
 \end{equation}
   where we denoted by $\vec E_N$ the vector
$(E_N, \ldots, E_N)$.

         We will also need a more rough estimate of the probability
             (\ref{mia})  out of
   the set ${\cal
       R}_{N,l}^{\eta}$.
       Let now the matrix $B_N(\o^{N,1},\ldots, \o^{N,l})$
   be of the rank $r\leq l$.
      Then, if $r<l$, there are $r$ paths among
$\o^{N,1},\ldots, \o^{N,l}$ such that
   corresponding $r$ random variables  $\eta(\o^{N,i})$ form the
   basis.
        Without loss of generality we may assume that these
   are $\o^{N,1},\ldots, \o^{N,r}$.
      Then the matrix
$B_N(\o^{N,1},\ldots, \o^{N,r})$
  is non-degenerate
  and $\eta(\o^{N,r+1}),\ldots,\eta(\o^{N,l})$
  are linear combinations of
$\eta(\o^{N,1}),\ldots,\eta(\o^{N,r})$.
  We may now estimate from above the probabilities (\ref{zet})
    by the probabilities
    $\P(\forall_{i=1}^{r} :  |\eta(\o^{N,i})-E_N|<b_i \delta_N)
$ that can be expressed in terms of the $r$-dimmensional
   integrals
like (\ref{mia}). Consequently, in this case
\begin{equation}
\P(\forall_{i=1}^{l} : |\eta(\o^{N,i})-E_N|<b_i \delta_N)
  \leq  \frac{ (2\delta_N/\sqrt{2 \pi})^r b_1\cdots b_r }
      {\sqrt{ {\rm det }B_N(\o^{N,1},\ldots, \o^{N,r}) }}.
       \end{equation}
   From the definition of the matrix elements,
      one sees that ${\rm det }B_N(\o^{N,1},\ldots, \o^{N,l})$
        is a finite polynomial in the variables $1/N$.
          Hence, if the rank of $B(\o^{N,1},\ldots, \o^{N,r})$
          equals $r$, we have for all $N>0$
\begin{equation}
\label{es2}
 \P(\forall_{i=1}^{l} : |\eta(\o^{N,i})-E_N|<b_i
\delta_N) \leq
     2^{-Nr}e^{c^2 r N^{2\alpha}/2}N^{k(r)}
 \end{equation}
 for some $k(r)>0$.

\medskip

 \noindent{\it Step 3.}
      Armed with (\ref{zgu}), (\ref{es1}) and (\ref{es2}),
    we now proceed with the proof of the theorem.

  For given $\alpha \in ]0, 1/4[$, let us choose
   first
    $\eta_0 \in ]0, 1/4[$ such that
\begin{equation}
\label{eta_0}
  2\alpha-1/2+\eta_0<0.
  \end{equation}
     Next, let us choose $\eta_1>\eta_0$
       such that
\begin{equation}
\label{keta_0}
  2\alpha-1/2+\eta_1<2\eta_0,
  \end{equation}
  then $\eta_2>\eta_1$ such that
\begin{equation}
\label{eta_1}
  2\alpha-1/2+\eta_2<2\eta_1,
  \end{equation}
etc. After $i-1$ steps we choose  $\eta_i >\eta_{i-1}$ such that
\begin{equation}
\label{eta_i}
  2\alpha-1/2+\eta_i<2\eta_{i-1}.
  \end{equation}
   Let us take e.g.\ $\eta_i=(i+1)\eta_0$.
    We stop the procedure at
  $n = [\alpha/\eta_0]$th step, that is
\begin{equation}
 \label{eta_n}
  n=\min\{i\geq 0 : \alpha <\eta_i\}.
\end{equation}
  Note that $\eta_{n-1}\leq \alpha<1/4$, and then
   $\eta_n=\eta_{n-1}+\eta_0<1/2$.

    We will prove that the sum
(\ref{zet}) over ${\cal R}_{N,l}^{\eta_0}$
  converges to $b_1\cdots b_l$, while those over
   ${\cal R}_{N,l}^{\eta_i}\setminus {\cal R}_{N,l}^{\eta_{i-1}}$
     for $i=1,2,\ldots,n$ and the one over
$S_N^{\otimes l} \setminus {\cal R}_{N,l}^{\eta_{n}}$
  converge o zero.

By (\ref{es1}), each term of the sum (\ref{zet})
           over ${\cal R}^{\eta_0}_{N,l}
 $ equals
$$(2\delta_N/\sqrt{2\pi})^l (b_1\cdots b_l)
e^{- \|\vec E_N\|^2 (1+O(N^{\eta_0-1/2}))/2 }(1+o(1)).
$$
   Here $e^{\|\vec E_N\|^2 \times O(N^{\eta_0-1/2})}
      =1+o(1)$ by the choice  (\ref{eta_0}) of $\eta_0$.
  Then, by the definition of $\delta_N$
    (\ref{delta}),  each term of the sum (\ref{zet})
           over ${\cal R}^{\eta_0}_{N,l}$ is
$$ (b_1\cdots b_l) 2^{-Nl}(1+o(1))$$
  uniformly for $(\omega^{N,1},\ldots, \o^{N,l}) \in {\cal
       R}_{N,l}^{\eta_0}$.
         The number of terms in this
           sum is $|{\cal
       R}_{N,l}^{\eta_0}|$, that is
         $2^{Nl}(1+o(1))$ by (\ref{zgu}).
     Hence, the sum (\ref{zet}) over
 ${\cal R}^{\eta_0}_{N,l}
 $ converges to $b_1\cdots b_l$.

    Let us  consider the sum over
${\cal R}_{N,l}^{\eta_i}\setminus {\cal R}_{N,l}^{\eta_{i-1}}$
     for $i=1,2,\ldots,n$.
        Each term in this sum equals
$$(2\delta_N/\sqrt{2\pi})^l (b_1\cdots b_l)
e^{- \|\vec E_N\|^2 (1+O(N^{\eta_i-1/2})/2 }(1+o(1))
$$
uniformly for $(\omega^{N,1},\ldots, \o^{N,l}) \in {\cal
       R}_{N,l}^{\eta_i}$.  Then, by the definition
         of  $\delta_N$ (\ref{delta}),  it is bounded by
        $2^{-Nl} C_i e^{h_i N^{2\alpha -1/2+\eta_i}}$
   with some constants $C_i, h_i>0$.
    The number of terms in this sum
is not greater than $|S_{N}^{\otimes l} \setminus {\cal
R}_{N,l}^{\eta_{i-1}}|$
   which is bounded due to  (\ref{zgu})
     by $C N 2^{Nl}\exp(-h N^{2\eta_{i-1}})$.
   Then by the choice of $\eta_i$
        (\ref{eta_i}) this sum converges to zero
     exponentially fast.

 Let us now treat the sum over
$S_N^{\otimes l} \setminus {\cal R}_{N,l}^{\eta_{n}}$.
     Let us first study the sum
over $(\o^{N,1},\ldots, \o^{N,l})$ such that
  the matrix $B_N(\o^{N,1},\ldots, \o^{N,l})$ is non-degenerate.
       By (\ref{es2})  each term in this sum
   is bounded by
 $ 2^{-Nl}e^{c^2 l N^{2\alpha}/2}N^{k(l)}$
   for some $k(l)>0$.
         The number of terms in this sum is bounded by
  $ C N 2^{Nl}\exp(-h N^{2\eta_{n}})$ by (\ref{zgu}). Since
     $\alpha<\eta_n$ by (\ref{eta_n}),
     this sum converges to zero exponentially fast.

     Let us finally turn to
the sum over $(\o^{N,1},\ldots, \o^{N,l})$ such that
  the matrix $B(\o^{N,1},\ldots, \o^{N,l})$
   is degenerate of the rank $r<l$.
By (\ref{es2})  each term in this sum
   is bounded by
   \begin{equation}
   \label{pp}
  2^{-Nr}e^{c^2 r N^{2\alpha}/2}N^{k(r)}
  \end{equation}
   for some $k(r)>0$.

      There are $r$ paths among
$\o^{N,1},\ldots, \o^{N,l}$ such that
   corresponding $\eta(\o^{N,i})$ form the basis.
        Without loss of generality we may assume that these
   are $\o^{N,1},\ldots, \o^{N,r}$.
     Note that $\o^{N,1},\ldots, \o^{N,r}$
  are such that it can not be for no one
     $m \in [0,\ldots, N]$ that
         $\o^1_m,\ldots, \o^r_m$ are all different.
       In fact, assume that
 $\o^1_m,\ldots, \o^r_m$ are all different. Then
     $\eta(m, \o^{1}_m),\ldots, \eta(m, \o^{r}_m)$
     are independent identically distributed random variables
         and $\eta(m, \o^{r+1}_m)=
           \mu_1 \eta(m, \o^{1}_m)+\cdots + \mu_r \eta(m,
           \o^{r}_m)$.
    If $\o^{r+1}_m$ is different from all $\o^1_m,\ldots, \o^r_m$,
         then $\eta(m, \o^{r+1}_m)$ is independent from
           all of $\eta(m, \o^{1}_m),\ldots,\eta(m,
           \o^{r}_m)$, then  the linear coefficients, being the
             covariances of $\eta(m, \o^{r+1}_m)$
              with  $\eta(m, \o^{1}_m),\ldots, \eta(m, \o^{r}_m)$,
    are $\mu_1=\cdots=\mu_r=0$.
     So, $\eta(\o^{N,r+1})$
         can not be a non-trivial linear combination
    of $\eta(\o^{N,1}),\ldots, \eta(\o^{N,r})$.
 If $\o^{r+1}_m$ equals one of $\o^1_m,\ldots, \o^r_m$,
   say $\o^{i}_m$, then again by computing the
  covariances of  $\eta(m, \o^{r+1}_m)$
              with  $\eta(m, \o^{1}_m),\ldots, \eta(m, \o^{r}_m)$,
    we get $\mu_i=1$, $\mu_j=0$
     for $j=1,\ldots, i-1,i+1,\ldots,r$.
            Consequently,
   $\eta(\o^{i}_k)=\eta(\o^{r+1}_k)$
     for all $k=1,\ldots, N$, so that
        $\o^{N,i}=\o^{N,r+1}$. But this is impossible
          since the sum  (\ref{zet})
               is taken over \underline{different\/}  paths
          $\o^{N,1},\ldots, \o^{N,l}$.
  Thus the sum is taken only over paths
$\o^{N,1},\ldots, \o^{N,r}$ where at each moment of time
  at least two of them are at the same place.

      The number of such sets of $r$ different
            paths is exponentially smaller than
      $2^{Nr}$ : there exists $p>0$ such that
  is does not exceed $2^{Nr}e^{-pN}$.
(In fact, consider $r$  independent simple
   random walks on ${\bf Z}$ that at a given moment of time
     occupy any $k<r$ different points of ${\bf Z}$.
   Then with probability not less than
   $(1/2)^r$, at the next moment of time, they
      occupy at least $k+1$ different points.
   Then with probability not less than $((1/2)^r)^r$
 at least once during $r$ next moments of time
   they will occupy $r$ different points.
  So, the number of sets of different $r$ paths that
         at each moment of time during $[0,N]$
   occupy at most $r-1$ different points
    is not greater than $2^{Nr}(1-(1/2^{r})^r)^{[N/r]}$.)

    Given any set of $r$ paths with
      $\eta(\o^{N,1}),\ldots, \eta(\o^{N,r})$ linearly independent,
            there is an $N$-independent
    number of possibilities to complete it by linear
    combinations $\eta(\o^{N,r+1}),\ldots \eta(\o^{N,l})$.
      To see this, first consider the equation
 $\lambda_1\eta(\o^{N,1})+\cdots+\lambda_r\eta(\o^{N,r})=0$
   with unknown $\lambda_1,\ldots, \lambda_r$.
     For any moment of time $m \in [0,N]$
       this means $\lambda_1 \eta(m, \o_m^{1})+\cdots +\lambda_r
       \eta(m, \o_m^r)=0$.
             If $\o_m^{i_1}=\o_m^{i_2}=\cdots \o_m^{i_k}$
               but $\o_m^j\ne \o_m^{i_1}$ for all $
     j\in \{1,\ldots,r\}\setminus\{i_1,\ldots, i_k\}$,
     then $\lambda_{i_1}+\cdots + \lambda_{i_k}=0$.
         Then for any
           $m \in [0, N]$
   the equation $\lambda_1 \eta(m, \o_m^{1})+\cdots +\lambda_r
       \eta(m, \o_m^r)=0$ splits into a certain number
         $n(m)$ ($1\leq n(m) \leq r$)
   equations of  type $\lambda_{i_1}+\cdots + \lambda_{i_k}=0$.
       Let us construct a matrix $A$ with $r$ columns
   and at least $N$ and at most $rN$ rows in the following way.
      For any $m>0$, according to given $\o_m^1,\ldots, \o_m^r$,
      let us add to A  $n(m)$ rows : each equation
          $\lambda_{i_1}+\cdots + \lambda_{i_k}=0$ gives
       a row with $1$ at places $i_1,\ldots, i_k$ and
          $0$ at all other places.
          Then the equation
           $\lambda_1\eta(\o^{N,1})+\cdots+\lambda_r\eta(\o^{N,i})=0$
             is equivalent $A \vec \lambda =\vec 0$
               with $\vec \lambda=(\lambda_1,\ldots, \lambda_r)$.
      Since this equation has only a trivial solution $\vec \lambda=0$,
    then the rank of $A$ equals $r$.
         The matrix $A$ contains at most $2^r$ different rows.
              There is less than  $(2^r)^r$ possibilities
    to choose $r$ linearly independent of them.
    Let $A^{r \times r}$ be an $r \times r$
        matrix consisting of  $r$ linearly independent rows of $A$.
        The fact that $\eta(\omega^{N,r+1})$ is
     a linear combination
      $\mu_1\eta(\o^{N,1})+\cdots+\mu_r\eta(\o^{N,r})=\eta(\o^{N,r+1})$
      can be written as  $A^{r \times r} \vec \mu =\vec b$
   where the vector $\vec b$ contains only  $1$
   and $0$ : if a given row $t$ of the matrix
     $A^{r \times r}$ corresponds to the $m$th step
       of the random walks and has $1$ at places
   $i_1,\ldots,i_k$ and $0$ elsewhere, then
     we put $b_t=1$ if $\o_m^{i_1}=\o_m^{r+1}$
      and $b_t=0$ if $\o_m^{i_1}\ne \o_m^{r+1}$.
   Thus, given
     $\o^{N,1},\ldots, \o^{N,r}$,
           there is an $N$ independent number
     of possibilities to write the system $A^{r \times r} \vec \mu =\vec b$
        with non degenerate matrix $A^{r \times r}$
      which determines uniquely linear coefficients
         $\mu_1,\ldots, \mu_r$ and consequently
   $\eta(\o^{N,r+1})$ and the path $\o^{N,r+1}$
     itself through these linear coefficients.
   Hence, there is not more possibilities to
     choose $\o^{N,r+1}$
   than the number of  non-degenerate matrices
    $A^{r \times r}$ multiplied by the number of vectors $\vec
    b$, that is roughly not more than $2^{r^2+r}$.

  These observations lead to the fact that the sum (\ref{zet})
     with the  covariance matrix $B_N(\o^{N,1},\ldots, \o^{N,l})$ of the
       rank  $r$ contains at most $(2^{r^2+r})^{l-r} 2^{Nr}e^{-p N}$
  different terms with some constant $p>0$.
    Then, taking into account the estimate (\ref{pp})
      of each term with $2\alpha<1$, we deduce that it converges to zero
     exponentially fast.
 This finishes the proof
     of  (\ref{th1e}).

     To show (\ref{th1b}), we have been already noticed
       that the sum of terms
       $\P(\forall_{i=1}^{2} : |\eta(\o^{N,i})-E_N|<b_i \delta_N)$
         over all pairs  of different paths $\o^{N,1}, \o^{N,2}$
           in $S_N^{\otimes l} \setminus
                {\cal R}_{N,l}^{\eta_0}$
          converges to zero exponentially fast.
   Then (\ref{th1b}) follows from the Borel-Cantelli lemma.

\medskip

\noindent{\bf Proof of Theorem \ref{th2}.}
    We  have again to verify the hypothesis of Theorem \ref{tc}
         for $V_{i,M}$ given by
  $\delta_N^{-1}
  |\eta(\o^{N,i})-E_N|$, i.e.\ we must show (\ref{zet}).

  For  $\beta \in ]0, 1[$  let us denote by
 $$ {\cal K}_{N,l}^{\beta}=\{ (\o^{N,1},\ldots, \o^{N,l}) :
    \ {\rm cov}(\eta(\o^{N,i}),\eta( \o^{N,j}))\leq N^{\beta-1},\
    \forall i,j=1,\ldots,l,\
      i\ne j\}.$$

\noindent{\it Step 1.}  In this step we estimate the capacity of the
complementary
     set to ${\cal K}_{N,l}^{\beta}$ in (\ref{kk}) and (\ref{kkk}).

   We have :  \begin{eqnarray}
   \lefteqn{ |S_N^{\otimes,l}\setminus {\cal K}_{N,l}^{\beta}|
 }\label{zgu1} \\
      &\leq & (l(l-1)/2) (2d)^{N(l-2)}
        \#\Big\{\omega^{N,1},\omega^{N,2} :
                  \#\{m \in[0,\ldots,N] :
                \omega_m^1 - \o_m^2=0\} > N^{\beta}\Big\}.
         \nonumber
\end{eqnarray}
  It has been shown in the proof of Theorem \ref{th1} that
   the number
$$\#\Big\{\omega^{N,1},\omega^{N,2} :
                  \#\{m \in[0,\ldots,N] :
                \omega_m^1 - \o_m^2=0\} > N^{\beta}\Big\}$$
 equals the number of paths of a simple random walk
 within the period $[0,2N]$ that visit the origin
  at least $[N^\beta]+1$ times.

      Let $W_r$ be the time of the $r$th return to the origin
  of a simple random walk
    ($W_1=0$), $R_N$ be the number of returns
   to the origin in the first $N$ steps.
          Then for any integer $q$
      $$P(R_N \leq q)=P(W_1+(W_2-W_1)+\cdots +(W_q-W_{q-1}) \geq N)
   \geq \sum_{k=1}^{q-1} P(E_k)$$
     where $E_k$ is the event that exactly $k$ of the variables
        $W_s-W_{s-1}$ are greater or equal than $N$,
          and $q-1-k$ are less than $N$. Then
$$\sum_{k=1}^{q-1} P(E_k)=\sum_{k=1}^{q-1} {q-1 \choose k}
     P(W_2-W_1 \geq N)^k (1-  P(W_2-W_1 \geq N))^{q-1-k}$$
     $$=
       1-(1-  P(W_2-W_1 \geq N))^{q-1}.$$
    It is shown in \cite{ET}
    that in the case $d=2$
    $$P(W_2-W_1 \geq N)
       =\pi (\log N)^{-1}(1+ O((\log N)^{-1})), \ \ \ N \to \infty.$$
  Then
  $$ P(R_N >q) \leq \Big(1-\pi (\log N)^{-1}(1+o(1))\Big)^{q-1}.$$
  Consequently,
  $$ \#\Big\{\omega^{N,1},\omega^{N,2} :
                  \#\{m \in[0,\ldots,N] :
                \omega_m^1 - \o_m^2=0\} > N^{\beta}\Big\}
                $$
                $$=(2d)^{2N} P(R_{2N}>[N^\beta])$$
                  $$\leq
                 (2d)^{2N} \Big(1-\pi (\log 2N)^{-1}(1+o(1)) \Big)^{[N^\beta]-1}
                   \leq (2d)^{2N} \exp(- h (\log 2N)^{-1} N^{\beta}) $$
      with some constant $h>0$.
         Finally for $d=2$ and all $N>0$
           by (\ref{zgu1})
 \begin{eqnarray}
    |S_N^{\otimes l}\setminus {\cal K}_{N,l}^{\eta}|
  \leq  (2d)^{lN} \exp(- h_2 (\log 2N)^{-1} N^{\beta}) \label{kk}
  \end{eqnarray}
   with some constant $h_2>0$.

   In the case $d\geq 3$ the random walk is transient and
 $$P(W_2-W_1 \geq N)\geq P(W_2-W_1 =\infty)=\gamma_d>0.$$
   It follows that $\P(R_N>q)\leq (1-\gamma_d)^{q-1}$ and
   consequently
\begin{eqnarray}
    |S_N^{\otimes,l}\setminus {\cal K}_{N,l}^{\beta}|
  \leq  (2d)^{lN} \exp(- h_d N^{\beta}) \label{kkk}
  \end{eqnarray}
  with some constant $h_d>0$.

\medskip

 \noindent{\it Step 2.}    Proceeding exactly as in the proof of Theorem
     \ref{th1}, we obtain that uniformly for
        $(\omega^{N,1},\ldots, \o^{N,l}) \in {\cal
       K}_{N,l}^{\beta}$,
     \begin{equation}
     \label{est1}
     \P(\forall_{i=1}^{l} : |\eta(\o^{N,i})-E_N|<b_i \delta_N)
        $$
        $$ =(2\delta_N/\sqrt{2\pi})^l (b_1\cdots b_l)
   e^{-\|E_N\|^2(1+O(N^{\beta-1}))/2}(1+o(1))
 \end{equation}
   where we denoted by $\vec E_N$ the vector
$(E_N, \ldots, E_N)$.
   Moreover, if
the covariance the matrix $B_N(\o^{N,1},\ldots, \o^{N,l})$
   is of the rank $r\leq l$
     (using the fact that its determinant is a finite polynomial
        in the variables $1/N$)
            we get as in the proof of Theorem \ref{th1} that
\begin{equation}
\label{est2}
 \P(\forall_{i=1}^{l} : |\eta(\o^{N,i})-E_N|<b_i
\delta_N) \leq
     (2d)^{-Nr}e^{c^2 r N^{2\alpha}/2}N^{k(r)}
 \end{equation}
 for some $k(r)>0$.

\medskip

 \noindent{\it Step 3.}
   Having (\ref{kk}), (\ref{kkk}), (\ref{est1}) and (\ref{est2}),
      we are able to carry out the proof of the theorem.
    For given $\alpha \in ]0, 1/2[$, let us choose
   first $\beta_0>0$ such that
\begin{equation}
\label{beta_0}
  2\alpha-1+\beta_0<0.
  \end{equation}
     Next, let us choose $\beta_1>\beta_0$
       such that
\begin{equation}
\label{bketa_0}
  2\alpha-1+\beta_1<\beta_0,
  \end{equation}
  then $\beta_2>\beta_1$ such that
\begin{equation}
\label{beta_1}
  2\alpha-1+\beta_2<\beta_1,
  \end{equation}
etc. After $i-1$ steps we choose  $\beta_i >\beta_{i-1}$ such that
\begin{equation}
\label{beta_i}
  2\alpha-1+\beta_i<\beta_{i-1}.
  \end{equation}
   Let us take e.g.\ $\beta_i=(i+1)\beta_0$.
    We stop the procedure at
  $n = [2\alpha/\beta_0]$th step, that is
\begin{equation}
 \label{beta_n}
  n=\min\{i\geq 0 : 2\alpha <\beta_i\}.
\end{equation}
  Note that $\beta_{n-1} \leq 2\alpha$, and then
   $\beta_n=\beta_{n-1}+\beta_0<2\alpha+ 1-2\alpha=1$.

    We will prove that the sum
(\ref{zet}) over ${\cal K}_{N,l}^{\beta_0}$
  converges to $b_1\cdots b_l$, while those over
   ${\cal K}_{N,l}^{\beta_i}\setminus {\cal K}_{N,l}^{\beta_{i-1}}$
     for $i=1,2,\ldots,n$ and the one over
$S_N^{\otimes l} \setminus {\cal K}_{N,l}^{\beta_{n}}$
  converge o zero.

By (\ref{est1}), each term of the sum (\ref{zet})
           over ${\cal K}^{\beta_0}_{N,l}
 $ equals
$$(2\delta_N/\sqrt{2\pi})^l (b_1\cdots b_l)
e^{- \|\vec E_N\|^2 (1+O(N^{\beta_0-1}))/2 }(1+o(1)).
$$
   Here $e^{\|\vec E_N\|^2 \times O(N^{\beta_0-1})}
      =1+o(1)$ by the choice  (\ref{beta_0}) of $\beta_0$.
  Then, by the definition of $\delta_N$
    (\ref{delta1}), each term of the sum (\ref{zet})
           over ${\cal K}^{\beta_0}_{N,l}$ is
$$ (b_1\cdots b_l) (2d)^{-Nl}(1+o(1))$$
  uniformly for $(\omega^{N,1},\ldots, \o^{N,l}) \in {\cal
       K}_{N,l}^{\eta_0}$.
         The number of terms in this
           sum is $|{\cal
       K}_{N,l}^{\beta_0}|$, that is
         $(2d)^{Nl}(1+o(1))$ by (\ref{kk}) and (\ref{kkk}).
     Hence, the sum (\ref{zet}) over
 ${\cal K}^{\beta_0}_{N,l}
 $ converges to $b_1\cdots b_l$.

    Let us  consider the sum over
${\cal K}_{N,l}^{\beta_i}\setminus {\cal K}_{N,l}^{\beta_{i-1}}$
     for $i=1,2,\ldots,n$. By (\ref{est1})
        each term in this sum equals
$$(2\delta_N/\sqrt{2\pi})^l (b_1\cdots b_l)
e^{- \|\vec E_N\|^2 (1+O(N^{\beta_i-1})/2 }(1+o(1))
$$
uniformly for $(\omega^{N,1},\ldots, \o^{N,l}) \in {\cal
       K}_{N,l}^{\beta_i}$.  Then, by the definition
         of  $\delta_N$ (\ref{delta1}),  it is bounded by
       the quantity $(2d)^{-Nl} C_i e^{h_i N^{2\alpha -1+\beta_i}}$
   with some constants $C_i, h_i>0$.
    The number of terms in this sum
is not greater than $|S_{N}^{\otimes l} \setminus {\cal
K}_{N,l}^{\beta_{i-1}}|$
   which is bounded
     by $(2d)^{Nl}\exp(-h_2 N^{\beta_{i-1}} (\log 2 N)^{-1})$
       in the case $d=2$ due to (\ref{kk})
        and
  by the quantity $(2d)^{Nl}\exp(-h_d N^{\beta_{i-1}} )$
  in the case $d\geq 3$ due to (\ref{kkk}).
   Then by the choice of $\beta_i$
        (\ref{beta_i}) this sum converges to zero
     exponentially fast.

 Let us now treat the sum over
$S_N^{\otimes l} \setminus {\cal K}_{N,l}^{\beta_{n}}$.
     Let us first analyze the sum
over $(\o^{N,1},\ldots, \o^{N,l})$ such that
  the matrix $B_N(\o^{N,1},\ldots, \o^{N,l})$ is non-degenerate.
       By (\ref{est2})  each term in this sum
   is bounded by
 $ (2d)^{-Nl}e^{c^2 l N^{2\alpha}/2}N^{k(l)}$
   for some $k(l)>0$.
         The number of terms in this sum is bounded by
 the quantity  $ (2d)^{Nl}\exp(-h_2 N^{\beta_{n}} (\log 2N)^{-1})$
  in the case $d=2$ and by
 $ (2d)^{Nl}\exp(-h_d N^{\beta_{n}})$
  in the case $d\geq 3$ respectively by (\ref{kk}) and (\ref{kkk}) . Since
     $2\alpha<\beta_n$ by (\ref{beta_n}),
     this sum converges to zero exponentially fast.

Let us finally turn to the sum over $(\o^{N,1},\ldots, \o^{N,l})$
such that
  the matrix $B_N(\o^{N,1},\ldots, \o^{N,l})$
   is degenerate of the rank $r<l$.
By (\ref{est2})  each term in this sum
   is bounded by $
  (2d)^{-Nr}e^{c^2 r N^{2\alpha}/2}N^{k(r)}$
   for some $k(r)>0$, while exactly by
    the same arguments as in the proof of Theorem \ref{th1},
     (they are, indeed, valid in all dimensions)
    the number of terms in this sum is
    less than $O((2d)^{Nr})e^{-p N}$
    with some constant $p>0$.
       Hence, this last sum converges to zero
         exponentially fast as $2\alpha <1$.
             This finishes the proof  of (\ref{th2e}).
          The proof of (\ref{th2b}) is completely
           analogous to the one of
(\ref{th1b}).

\medskip

\noindent{\bf Proof of Theorem \ref{th0}.}
     We again concentrate on  the proof in the sum (\ref{zet})
       with $E_N=c$.

 \noindent{\it Step 1.}  First of all,  we  need a rather rough estimate
    of the probabilities of (\ref{zet}).
  Let $(\o^{N,1},\ldots, \o^{N,r})$
    be such that the matrix $B_N(\o^{N,1},\ldots, \o^{N,r})
     $ is  non-degenerate.
       We  prove in this step that there exists  a constant $k(r)>0$
       such that for any $N>0$ and any $(\o^{N,1},\ldots, \o^{N,r})$
    with non-degenerate $B_N(\o^{N,1},\ldots, \o^{N,r})$, we have:
\begin{equation}
\label{gg2} \P(\forall_{i=1}^{r} :  |\eta(\o^{N,i})-c|<b_i \delta_N)
   \leq (2d)^{-Nr}N^{k(r)}.
\end{equation}
   Let
   $$f^{\o^{N,1},\ldots, \o^{N,r}}_N(t_1,\ldots, t_r)
      = \E \exp\Big( i\sum_{k=1}^r t_k \eta(\o^{N,k})
      \Big)$$
 be the Fourier transform  of  $(\eta(\o^{N,1}),\ldots,
 \eta(\o^{N,r}))$. Then
\begin{eqnarray}
 \lefteqn{\P(\forall_{i=1}^{r} :  |\eta(\o^{N,i})-c|<b_i \delta_N)}
 \nonumber\\
&= &\frac{1}{(2\pi)^r}
  \int\limits_{{\bf R}^r} f^{\o^{N,1},\ldots, \o^{N,r}}_N(\vec t)
      \prod_{k=1}^r \frac{ e^{-i t_k (-b_k \delta_N+c)} -  e^{-i t_k (b_k
      \delta_N +c)}}{it_k} dt_1\cdots d t_r \label{fou}
\end{eqnarray}
   provided that the integrand is in $L^1({\bf R}^d)$.
  We will show that this is the case  due to the assumption
   made on $\phi$  and deduce the bound (\ref{gg2}).

     We know that the function $f_N^{\o^{N,1},\ldots, \o^{N,r}}(\vec
  t)$ is the product of $N$ generating functions :
   \begin{equation}
   \label{zey}
      f^{\o^{N,1},\ldots, \o^{N,r}}_N(\vec t)=
      \prod_{n=1}^N \E \exp\Big( i N^{-1/2}
         \sum_{k=1}^r t_k  \eta(n, \o^{N,k}_n)
      \Big).
      \end{equation}
   Moreover, each of these functions
   is itself a product of (at minimum $1$ and at maximum $r$)
     generating functions of type $\phi((t_{i_1}+\cdots
     +t_{i_k})N^{-1/2})$.
   More precisely, let us construct the matrix $A$ with
      $r$ columns and at least $N$ and at most $rN$
   rows as in the proof of Theorem~\ref{th1}. Namely,
          for each step $n=0,1,2,\ldots,N$,
    we add to the matrix $A$ at least $1$ and at most $r$ rows
  according to the following rule:
      if $\o_n^{N,i_1}=\o_n^{N,i_2}=\cdots=\o_n^{N,i_k}$
          and $\o_n^{N,j}\ne \o_n^{N,i_1}$
            for any $j\in \{1,\ldots,r\}
            \setminus \{i_1,\ldots,i_k\}$, we add to $A$
              a row with $1$ at places $i_1,\ldots, i_k$ and
       $0$ at other $r-k$ places.
  Then
   \begin{equation}
   \label{aa}
f^{\o^{N,1},\ldots, \o^{N,r}}_N(\vec t)=
  \prod_{j} \phi (N^{-1/2}(A \vec t)_j).
   \end{equation}
  Since $B_N(\o^{N,1},\ldots, \o^{N,r})
     $ is non-degenerate, the rank of the matrix $A$
       equals $r$.  Let us choose in $A$
         any $r$ linearly independent rows,
           and let us denote by $A^r$
      the $r \times r$ matrix
        constructed by them.
  Then by the assumption made on $\phi$
\begin{equation}
   \label{aaa}
|f^{\o^{N,1},\ldots, \o^{N,r}}_N(\vec t)| \leq
  \prod_{j=1}^r |\phi (N^{-1/2}(A^r \vec t)_j)|
    \leq \prod_{j=1}^r  \min \Big(1, \frac{ C N^{1/2} }{ |(A^r \vec
    t)_j|} \Big) \leq C^r N^{r/2}
       \prod_{j=1}^r  \min \Big(1, \frac{1 }{ |(A^r \vec
    t)_j|} \Big)
   \end{equation}
 with some constant $C>0$. Furthermore
\begin{equation}
\label{bb}
 \Big|\prod_{k=1}^r \frac{ e^{-i t_k (-b_k \delta_N+c)} -  e^{-i t_k (b_k
      \delta_N +c)}}{it_k} \Big| \leq
         \prod_{k=1}^r \min \Big( (2\delta_N)b_k, \
           \frac{2}{|t_k|}\Big) \leq C'\prod_{k=1}^r
  \min \Big((2d)^{-N}, \frac{1}{|t_k|} \Big)
\end{equation}
  with some $C'>0$. Hence,
  \begin{eqnarray}
 \lefteqn{\frac{1}{(2\pi)^r}
  \int\limits_{{\bf R}^r} \Big|f^{\o^{N,1},\ldots, \o^{N,r}}_N(\vec t)
      \prod_{k=1}^r \frac{ e^{-i t_k (-b_k \delta_N+c)} -  e^{-i t_k (b_k
      \delta_N +c)}}{it_k}\Big| dt_1\cdots d t_r }\nonumber\\
 & \leq & C_0 N^{r/2}
       \int
\prod_{k=1}^r
  \min \Big((2d)^{-N}, \frac{1}{|t_k|} \Big)
     \min \Big(1, \frac{1 }{ |(A^r \vec
    t)_k|} \Big) d \vec t \label{ss}
    \end{eqnarray}
  with some constant $C_0>0$  depending on the function $\phi$ and on
    $b_1,\ldots, b_r$ only.
       Since the matrix $A^r$ is non-degenerate,
     using easy arguments
    of linear algebra, one can show that for some
      constant $C_1>0$ depending on the matrix $A^{r}$ only,
  we have
\begin{equation}
 \int \prod_{k=1}^r
  \min \Big((2d)^{-N}, \frac{1}{|t_k|} \Big)
     \min \Big(1, \frac{1 }{ |(A^r \vec
    t)_k|} \Big)d\vec t
 \leq   C_1 \int \prod_{k=1}^r
  \min \Big((2d)^{-N}, \frac{1}{|t_k|} \Big)
     \Big(1, \frac{1 }{ |
    t_k|} \Big) d\vec t. \label{sdg}
\end{equation}
   The proof of (\ref{sdg}) is given in Appendix.
But the right-hand of (\ref{sdg}) is finite.
  This shows that the integrand in (\ref{fou})
     is in $L^1({\bf R}^d)$ and
   the inversion formula (\ref{fou}) is valid.
  Moreover, the right-hand side of (\ref{sdg})  equals $C_1 (2((2d)^{-N}
+  (2d)^{-N} N \ln 2d + (2d)^{-N}))^r$.
 Hence, the probabilities above are bounded by the quantity
   $C_0 N^{r/2} C_1 2^r(2+N \ln (2d))^r (2d)^{-Nr}$
     with $C_0$ depending on $\phi$ and $b_1,\ldots, b_r$ and
          $C_1$ depending on the choice of $A^r$.
To conclude the proof of (\ref{gg2}), it remains to remark
 that there is an $N$-independent number of possibilities
   to construct a matrix $A^{r}$ (at most $2^{r^2}$),
since it contains only $0$ or $1$.

\medskip

 \noindent{\it Step 2.}  We keep the notation $\RR_{N,l}^{\eta}$
   from (\ref{rr}) for $\eta \in ]0,1/2[$.
  The capacity of this set for $d=1$ is estimated in (\ref{zgu}).
    Moreover by (\ref{kk}) for $d=2$
 $$ |S_{N}^{\otimes l}\setminus \RR_{N,l}^{\eta}|=
  |S_{N}^{\otimes l}\setminus \KK_{N,l}^{\eta+1/2}|
     \leq (2d)^{Nl} \exp(-h_2 (\log 2N)^{-1} N^{1/2+\eta})$$
     and by (\ref{kkk}) for $d\geq 3$
  $$ |S_{N}^{\otimes l}\setminus \RR_{N,l}^{\eta}|
     =    |S_{N}^{\otimes l}\setminus \KK_{N,l}^{\eta+1/2}|
           \leq (2d)^{Nl}\exp(-h_d N^{1/2+\eta}),$$
  so that, for all $d\geq 1$ there are $h_d, C_d>0$
  such that for all $N> 0$
\begin{equation}
\label{rrs}
 |S_{N}^{\otimes l}\setminus \RR_{N,l}^{\eta}|
              \leq (2d)^{Nl}C_d N \exp(-h_d N^{2\eta}).
              \end{equation}

\medskip

\noindent{\it Sep 3.}   In this step we show that uniformly for
$(\o^{N,1},
     \ldots, \o^{N,l}) \in
\RR_{N,l}^{\eta}$
  \begin{equation}
 \label{gg}
\P(\forall_{i=1}^{l} :  |\eta(\o^{N,i})-c|<b_i \delta_N)=
(2d)^{-Nl}b_1\cdots b_l(1+o(1)).
\end{equation}

  For any $(\o^{N,1},\ldots \o^{N,l}) \in
\RR_{N,l}^{\eta}$, we can represent the probabilities in the
 sum (\ref{zet}) as
sums of four terms :
\begin{eqnarray}
 \lefteqn{\P(\forall_{i=1}^{l} :  |\eta(\o^{N,i})-c|<b_i \delta_N)}
 \nonumber\\
&= &\frac{1}{(2\pi)^l}
  \int_{{\bf R}^l} f^{\o^{N,1},\ldots, \o^{N,l}}_N(\vec t)
      \prod_{k=1}^l \frac{ e^{-i t_k (-b_k \delta_N+c)} -  e^{-i t_k (b_k
      \delta_N +c)}}{it_k} dt_1\cdots d t_l \nonumber\\
      & = &\sum_{m=1}^4 I_N^m(\o^{N,1},\ldots, \o^{N,l})\label{zuzu}
\end{eqnarray}
 where
 \begin{eqnarray}
  I_N^1& =& \frac{1}{(2\pi)^l} \int\limits_{{\bf R}^l}
       \prod_{k=1}^l \frac{ e^{-i t_k (-b_k \delta_N+c)} -  e^{-i t_k (b_k
      \delta_N +c)}}{it_k} e^{-\vec t B_N(\o^{N,1},\ldots, \o^{N,l})
      \vec t/2}d\vec t \label{dd} \\
       &&{}-\frac{1}{(2\pi)^l}  \int\limits_{\|t\|>\epsilon N^{1/6}}
       \prod_{k=1}^l \frac{ e^{-i t_k (-b_k \delta_N+c)} -  e^{-i t_k (b_k
      \delta_N +c)}}{it_k} e^{-\vec t B_N(\o^{N,1},\ldots, \o^{N,l})
      \vec t/2}d\vec t.\nonumber\\
  I_N^2& =& \frac{1}{(2\pi)^l} \int\limits_{\|t\|<\epsilon N^{1/6}}
   \prod_{k=1}^l \frac{ e^{-i t_k (-b_k \delta_N+c)} -  e^{-i t_k (b_k
      \delta_N +c)}}{it_k} \nonumber \\
       && \ \ \ \ \ \ \ \  \  \  \  \  \  \  \  \  \ \ \  \ \ {} \times  \Big(
    f^{\o^{N,1},\ldots, \o^{N,l}}_N(\vec t)-
       e^{-\vec t B_N(\o^{N,1},\ldots, \o^{N,l})
      \vec t/2}\Big) d\vec t\label{ff}\\
  I_N^3& =&\frac{1}{(2\pi)^l} \int\limits_{\epsilon N^{1/6}<\|t\|<\delta N^{1/2}}
   \prod_{k=1}^l \frac{ e^{-i t_k (-b_k \delta_N+c)} -  e^{-i t_k (b_k
      \delta_N +c)}}{it_k}
    f^{\o^{N,1},\ldots, \o^{N,l}}_N(\vec t)
     d\vec t\nonumber\\
  I_N^4& =& \frac{1}{(2\pi)^l}
  \int\limits_{ \|t\|>\delta N^{1/2}}
   \prod_{k=1}^l \frac{ e^{-i t_k (-b_k \delta_N+c)} -  e^{-i t_k (b_k
      \delta_N +c)}}{it_k}
    f^{\o^{N,1},\ldots, \o^{N,l}}_N(\vec t)
     d\vec t\nonumber
\end{eqnarray}
 with $\epsilon, \delta>0$ chosen according to the following
   Proposition \ref{pr1}.
\begin{pro}
\label{pr1}
 There exist constants $N_0,C,\epsilon, \delta, \zeta>0$ such
that
   for all $(\o^{N,1},\ldots \o^{N,l}) \in
\RR_{N,l}^{\eta}$  and all $N \geq N_0$
  the following estimates hold:
\begin{equation}
\label{base1} \Big|f^{\o^{N,1},\ldots, \o^{N,l}}_N(\vec t) -
e^{-\vec t B_N(\o^{N,1},\ldots, \o^{N,l}) \vec t/2} \Big|
  \leq \frac{C \|t\|^3}{\sqrt{N}}
  e^{ -\vec t B_N(\o^{N,1},\ldots, \o^{N,l}) \vec t/2},\ \ \ \
  \hbox{for all  } \|t\|\leq \epsilon N^{1/6}.
\end{equation}

\begin{equation}
\label{base2} \Big| f^{\o^{N,1},\ldots, \o^{N,l}}_N(\vec t)\Big|
\leq e^{-\zeta \|t\|^2} \ \ \
 \hbox{for all  } \|t\|<\delta \sqrt{N}.
\end{equation}
\end{pro}

     The proof of this proposition mimics the one of the Berry-Essen
  inequality and is given in Appendix.

    The first part of
$I_N^1$ is just the probability that
    $l$ Gaussian random variables with
        zero mean and covariance matrix
    $B_N(\o^{N,1},\ldots, \o^{N,l})$ belong
          to the intervals
 $[-\delta_N b_k+c, \delta_N b_k+c]$ for $k=1,\ldots, l$ respectively.
     This is
      \begin{equation}
    \label{mia0}
   \int\limits_{|z_j- c|\leq
   \delta_N b_j, \forall_{j=1}^l }
   \frac{e^{-(\vec z B^{-1}(\o^{N,1},\ldots,
 \o^{N,l})\vec z)/2}}{(2\pi)^{l/2}
   \sqrt{{\rm det} B(o^{N,1},\ldots,
 \o^{N,l})} }\, d\vec z
 $$
$$ =  (2\delta_N/\sqrt{2\pi})^l (b_1\cdots b_l) e^{-(\vec c
B^{-1}(\o^{N,1},\ldots,
 \o^{N,l})\vec c)/2}(1+o(1))$$
$$ =(2\delta_N/\sqrt{2\pi})^l (b_1\cdots b_l)
   e^{-lc^2(1+O(N^{\eta-1/2}))/2}(1+o(1))
   = (2d)^{-Nl}b_1\cdots b_l(1+o(1))
 \end{equation}
     uniformly for
        $(\omega^{N,1},\ldots, \o^{N,l}) \in {\cal
       R}_{N,l}^{\eta}$, where we denoted by $\vec c$ the vector
$(c, \ldots, c)$.
     Since
     \begin{equation}
     \label{pr}
   \prod_{k=1}^l \Big| \frac{ e^{-i t_k (-b_k \delta_N+c)} -  e^{-i t_k (b_k
      \delta_N +c)}}{it_k} \Big|
        \leq  (2 \delta_N b_1)\cdots (2\delta_N b_l)= O((2d)^{-Nl})
        \end{equation}
and the elements of the matrix $B_N(\o^{N,1},\ldots, \o^{N,l})$
out
of the
  diagonal are $O(N^{\eta-1/2})=o(1)$ as  $N \to \infty$,
     the second part of $I_N^1$
       is   smaller than
  $(2d)^{-Nl}$ exponentially (with exponential
    term $\exp(-h N^{1/3})$ for some $h>0$).

     There is a constant $C>0$ such that
    the term $I_N^2$ is bounded by
     $C (2d)^{-Nl} N^{-1/2}$
      for any $(\o^{N,1},\ldots \o^{N,l}) \in
\RR_{N,l}^{\eta}$ and all $N$ large enough.
     This follows from  (\ref{pr}),
       the estimate (\ref{base1})
     and again the fact that
   the elements of the matrix $B_N(\o^{N,1},\ldots, \o^{N,l})$ out of the
  diagonal are $O(N^{\eta-1/2})=o(1)$ as  $N \to \infty$.

     The third term $I_N^3$ is exponentially smaller than
     $(2d)^{-Nl}$
     by (\ref{pr}) and  the estimate (\ref{base2}).

    Finally, by (\ref{pr})
$$ |I_N^4|\leq (2 \delta_N b_1)\cdots (2\delta_N b_l)
        \int\limits_{\|t\|>\delta \sqrt{N}}
            |f^{\o^{N,1},\ldots, \o^{N,l}}_N(\vec t)| d\vec t=
              O((2d)^{-Nl})
             \int\limits_{\|t\|>\delta \sqrt{N}}
            |f^{\o^{N,1},\ldots, \o^{N,l}}_N(\vec t)| d\vec t.$$
  The function
   $ f^{\o^{N,1},\ldots, \o^{N,l}}_N(\vec t)$
   is the product of $N$ generating  functions (\ref{zey}).
         Note that for any pair $\o^{N,i}, \o^{N,j}$
  of  $(\o^{N,1},\ldots, \o^{N,l})\in \RR_{N,l}^{\eta}$,
              there are at most $N^{\eta+1/2}$
 steps $n$ where
            $\o^{N,i}_n= \o^{N,j}_n$.
         Then there are at least $N -[l(l-1)/2]N^{\eta+1/2}=a(N)$
          steps where all  $l$ coordinates  $\o^{N,i}$, $i=1,\ldots,
          l$,
               of the vector
$(\o^{N,1},\ldots, \o^{N,l}) \in \RR_{N,l}^{\eta}$
               are
           different.
           In this case
$$\E \exp\Big( i N^{-1/2}\sum_{k=1}^l t_k  \eta(n, \o^{N,k}_n)
\Big)=\phi(t_1 N^{-1/2})\cdots \phi(t_k N^{-1/2}).$$
    By the assumption made on $\phi$,
      this function is aperiodic and thus $|\phi(t)|<1$
        for $t\ne 0$.
    Moreover, for any $\delta>0$ there exists
      $h(\delta)>0$ such that $|\phi(t)|\leq 1-h(\delta)$
         for $|t|>\delta/l$.
  Then
  $$\int\limits_{\|t\|>\delta \sqrt{N}}
            |f^{\o^{N,1},\ldots, \o^{N,l}}_N(\vec t)| d\vec t
  \leq \int\limits_{\|t\| >\delta \sqrt{N}} |
    \phi(t_1 N^{-1/2})\cdots \phi(t_k N^{-1/2})|^{a(N)}d\vec t
$$
$$ = N^{l/2} \int\limits_{\|s\| >\delta  }
    |\phi(s_1 )\cdots \phi(s_k )|^{a(N)}d\vec s
      \leq N^{l/2}(1-h(\delta))^{a(N)-2}
     \int\limits_{\|s\| >\delta  }
    |\phi(s_1 )\cdots \phi(s_k )|^2d\vec s $$
      where $a(N)=N(1+o(1))$ and
         the last integral converges due to the assumption
       made on $\phi(s)$.
    Hence $I_N^4$ is exponentially smaller than $(2d)^{-Nl}$.
    This finishes the proof of (\ref{gg}).

\medskip

 \noindent{\it Step 4.} We are now able to prove the theorem using
  the estimates (\ref{gg2}),(\ref{rrs}) and  (\ref{gg}).
       By (\ref{gg}),
   the sum (\ref{zet}) over $\RR_{N,l}^{\eta}$
   (with fixed $\eta \in ]0,1/2[$) that contains
   by (\ref{rrs})$(2d)^{Nl}(1+o(1))$ terms,
    converges to  $b_1 \cdots b_l$.
     The sum (\ref{zet}) over $(\o^{N,1},\ldots, \o^{N,l}) \not \in
   \RR_{N,l}^{\eta}$ but with $B_N(\o^{N,1},\ldots, \o^{N,l})$
     non-degenerate, by (\ref{rrs}) has only at most
         $(2d)^{Nl} C N \exp(-h N^{2\eta})$ terms,
      while each of its terms by (\ref{gg2}) with $r=l$
   is of the order $(2d)^{-Nl}$ up to a polynomial term. Hence,
      this sum converges to zero.
   Finally,
due to the fact that in any
         set $(\o^{N,1},\ldots, \o^{N,l})$
    taken into account in (\ref{zet}) the paths
         are all different,
    the sum over $(\o^{N,1},\ldots, \o^{N,l}) \not \in
   \RR_{N,l}^{\eta}$ with
     $B_N(\o^{N,1},\ldots, \o^{N,l})$ of the rank $r<l$
      has an
     exponentially smaller number of terms than $(2d)^{Nr}$.
        This has been shown  in detail
        in the proof of Theorem \ref{th1} where
          the arguments did not depend on the dimension
  of the random walk.
        Since  by (\ref{gg2}) each of these terms
   is of the order $(2d)^{-Nr}$ up to a polynomial term,
     this sum converges to zero.
    This concludes the proof of (\ref{th1e0}).
      The proof of (\ref{th1b0}) is completely analogous
        to the one of (\ref{th1b}).

\section{Appendix}

\noindent{\it Proof of (\ref{sdg}).}
   It is carried out via trivial arguments of linear
   algebra.

   Let  $m=1,2,\ldots, r+1$,
       $D_{m-1}$ be a non-degenerate $r\times r$
   matrix  with the first $m-1$ rows having $1$ on the
  diagonal and $0$ outside of the diagonal.
    (Clearly, $D_0$ is just a non-degenerate matrix
      and  $D_r$ is the diagonal matrix with $1$ everywhere
     on the diagonal.)
  Let us introduce the integral
\begin{eqnarray}
 \lefteqn{ J^{m-1}(D_{m-1})}\nonumber\\
 & =& \int \prod_{k=1}^r
   \min \Big((2d)^{-N}, \frac{1}{|t_k|} \Big)
     \min \Big(1, \frac{1 }{ |(D_{m-1} \vec
    t)_k|} \Big)  d\vec t \nonumber \\
&=&
 \int \prod_{k=1}^{m-1}\min\Big((2d)^{-N}, \frac{1}{|t_k|} \Big)
          \min\Big(1, \frac{1}{|t_k|} \Big)
     \prod_{k=m}^r
  \min \Big((2d)^{-N}, \frac{1}{|t_k|} \Big)
     \min \Big(1, \frac{1 }{ |(D_{m-1} \vec
     t)_k|} \Big)  d\vec t.\nonumber
 \end{eqnarray}
  Sice $D_{m-1}$ is non-degenerate, there exists $i \in \{m,\ldots,r\}$
     such that $d_{m,i} \ne 0$
             and the matrix $D_m$
 which is obtained from the matrix $D_{m-1}$ by replacing its $m$th
 row
   by the one with $1$ at the place $(m,i)$ and $0$
     at all places $(m,j)$ for $j \ne i$
       is non-degenerate.
  Without loss of generality we may assume that $i=m$
    (otherwise juste permute the $m$th with the $i$th column
      in $D_{m-1}$ and $t_i$ with $t_m$ in the integral $J^{m-1}(D_{m-1})$ above).
 Since either
  $|t_{m-1}| < |(D_{m-1} \vec t)_{m-1}|$
   or $|t_{m-1}| \geq |(D_{m-1} \vec t)_{m-1}|$,
  we can estimate $J^{m-1}(D_{m-1})$ roughly by the sum of the following
   two terms :
\begin{eqnarray}
  \lefteqn{J^{m-1}(D_{m-1})}\nonumber \\
& \leq &
  \int \prod_{k=1}^{m}\min\Big((2d)^{-N}, \frac{1}{|t_k|}
\Big)
          \min\Big(1, \frac{1}{|t_k|} \Big)
     \prod_{k=m+1}^r
  \min \Big((2d)^{-N}, \frac{1}{|t_k|} \Big)
     \min \Big(1, \frac{1 }{ |(D_{m-1} \vec
    t)_k|} \Big)  d\vec t \nonumber\\
 &&{} +
  \int \prod_{k=1}^{m-1}\min\Big((2d)^{-N}, \frac{1}{|t_k|} \Big)
          \min\Big(1, \frac{1}{|t_k|} \Big)
\min\Big((2d)^{-N}, \frac{1}{|(D_{m-1} \vec t)_m|} \Big)
          \min\Big(1, \frac{1}{| (D_{m-1} \vec t)_m|} \Big)\nonumber \\
   &&\ \ \ \ \ \ \ \ {}\times  \prod_{k=m+1}^r
  \min \Big((2d)^{-N}, \frac{1}{|t_k|} \Big)
     \min \Big(1, \frac{1 }{ |(D_{m-1} \vec
    t)_k|} \Big)  d\vec t. \label{dd}
\end{eqnarray}
      The first term here is just $J^{m}( D_m)$.
  Let us make a change of variables
   in the second one : let $\vec z =B_{D_{m-1}} \vec t$,
    where the matrix $B_{D_{m-1}}$ is chosen such that
              $z_1=t_1,\ldots, z_{m-1}=t_{m-1},
          z_m =(D_{m-1} \vec t)_m, z_{m+1}=t_{m+1},
       \ldots, z_{r}=t_r$. (Clearly, its $m$th row
         is the same as in the matrix $D_{m-1}$,
           and it has $1$ on the diagonal in all other
               $r-1$ rows and $0$ outside of it.)
  Since $d_{m,m}\ne 0$, the matrix $B$ is non-degenerate.
      Then $ D_{m-1} \vec t=D_{m-1} B^{-1}_{D_{m-1}} \vec z$,
  where the matrix $D_{m-1} B^{-1}_{D_{m-1}} $ is non-degenerate,
   and, moreover, it has the first $m$ rows with $1$ on the
   diagonal and $0$ outside of it, as we have
      $(D_{m-1} \vec t)_1=t_1=z_1, \ldots, (D_{m-1} \vec t)_{m-1}=t_{m-1}
        =z_{m-1}, (D_{m-1} \vec t)_m=z_m$.
  Then (\ref{dd}) can be written as
\begin{equation}
\label{ddd}
 J^{m-1}(D_{m-1}) \leq J^m(D_{m})+ d_{m,m}^{-1}
    J^m( D_{m-1} B^{-1}_{D_{m-1}}).
\end{equation}

   Now, observe that the left-hand side  of (\ref{sdg}) is
        $J^0(A^r)$.
  By (\ref{ddd}) it is bounded by
   $J^1(A^r_1)+ a_{1,1}^{-1} J^1(A^r B^{-1}_{A^r})$.
    Again by (\ref{ddd}) each of these two terms can be estimated
        by a sum of two terms
  of type $J^2(\cdot)$ etc. After $2^r$ applications
    of (\ref{ddd}) $J^{0}(A^r)$
    is bounded by a sum of $2^r$ terms
         of type $J^r(D_r)$ multiplied by some constants
    depending only on the initial matrix $A_r$.
        But all these $2^r$ terms $J^r(D_r)$ are the same
   as in the right-hand side of (\ref{sdg}).

\medskip

\noindent{\it Proof of Proposition \ref{pr1}.}
   We use the representation (\ref{aa}) of
     $f_N^{\o^{N,1},\ldots, \o^{N,l}}(\vec t)
   $ as the product of a certain number
      $K(N,\o^{N,1},\ldots, \o^{N,l})$
         (denote it shortly by $K(N,\o)$, clearly
          $ N \leq K(N, \o) \leq lN$)
             of generating functions $ \phi(N^{-1/2}(A \vec t)_j)$
    where at most $2^l$ are different.
      Each of them is
        of the form $\E \exp(i N^{-1/2}(t_{i_1}+\cdots +t_{i,k})X)$
     with $X$ a standard Gaussian random variable.
  Applying the fact that $|e^{iz}-1-iz-(iz)^2/2!|\leq |z|^3/3!$
    for any $z\in {\bf R}$, we can write
\begin{equation}
\label{jj} \phi(N^{-1/2}(A \vec t)_j)
  = 1- \frac{ ((A \vec t)_j)^2}{2! N}-\theta_j \frac{((A \vec t)_j)^3
      \E|X|^3}{ 3!
  N^{3/ 2}} \equiv 1-\alpha_j
  \end{equation}
     with some complex $\theta_j$ with $|\theta_j|<1$.
      It follows that there are some constants $C_1,C_2>0$ such that
         for any $( \o^{N,1},\ldots, \o^{N,l}  ) \in
         \RR_{N,l}^{\eta}$ and any $j$
 we have:       $|\alpha_j|\leq C_1 \|\vec t\|^2 N^{-1}+ C_2\|\vec t\|^3
       N^{-3/2}$. Then $|\alpha_j|<1/2$ and $|\alpha_j|^2 \leq
         C_3 \|\vec t\|^3
       N^{-3/2}$  with some $C_3>0$
                 for all $\vec t$ of the absolute
              value $\|\vec t\|\leq \delta \sqrt{N}$
           with $\delta>0$ small enough.
            Thus $ \ln \phi(N^{-1/2}(A \vec t)_j)
              =-\alpha_j+\tilde \theta_j \alpha_j^2/2$
  (using the expansion $\ln(1+z)=z +\tilde \theta z^2/2$
   with some $\tilde \theta$
     of the absolute value $|\tilde \theta|<1$
        which is true for all $z$ with $|z|<1/2$)
        for all $( \o^{N,1},\ldots, \o^{N,l}  ) \in
         \RR_{N,l}^{\eta}$ and for all $\vec t$ with
          $\|\vec t\|\leq \delta \sqrt{N}$
        with some $\tilde \theta_j$ such that $|\tilde \theta_j|<1$.
          It follows that
\begin{equation}
\label{zgh}
      f_N^{\o^{N,1},\ldots, \o^{N,l}}(\vec t)
   = \exp\Big( -\sum_{j=1}^{K(N,\o)}
        \alpha_j+ \sum_{j=1}^{K(N,\o)}
            \tilde \theta_j \alpha_j^2/2\Big).
            \end{equation}
          Since $A^*A=B_N(\o^{N,1},\ldots,
       \o^{N,l})$,   here  $-\sum_{j=1}^{K(N,\o)}
        \alpha_j =-\vec t B_N(\o^{N,1},\ldots, \o^{N,l}) \vec t/2
           +\sum_{j=1}^{K(N,\o)} p_j$ where
             $|p_j| \leq C_2\|\vec t\|^3
       N^{-3/2}$.  Then
\begin{equation}
\label{zgh1}
    f_N^{\o^{N,1},\ldots, \o^{N,l}}(\vec t)=
       \exp \Big( -\vec t B_N(\o^{N,1},\ldots, \o^{N,l}) \vec t/2 \Big)
          \exp\Big(\sum_{j=1}^{K(N,\o)}p_j+
            \tilde \theta_j \alpha_j^2/2\Big)
\end{equation}
              where
           $|p_j|+  |\tilde \theta_j \alpha_j^2/2| \leq
      (C_2+C_3/2)\|\vec t\|^3 N^{-3/2}$ for all $j$.
   Since $K(N,\o) \leq l N$, we have
   \begin{equation}
   \label{qk}
     \Big| \sum_{j=1}^{K(N,\o)}
         p_j+  \tilde \theta_j \alpha_j^2/2 \Big|\leq (C_2+C_3/2)l\|t\|^3
      N^{-1/2}.
      \end{equation}
          It follows that for $\epsilon>0$ small enough
         $| \exp (\sum_{j=1}^{K(N,\o)}
         p_j+  \tilde \theta_j \alpha_j^2/2)-1|
           \leq C_4 \|\vec t\|^3 N^{-1/2}$
             for all $\vec t$ with $\|\vec t\| \leq \epsilon N^{1/6}$.
          This proves (\ref{base1}).
   Finally
\begin{equation}
\label{zlt}
 |f_N^{\o^{N,1},\ldots, \o^{N,l}}(\vec t)| \leq
       \exp \Big( -\vec t B_N(\o^{N,1},\ldots, \o^{N,l}) \vec t/2 \Big)
    \exp \Big((C_2+C_3/2)l\|t\|^3
      N^{-1/2}\Big).
   \end{equation}
Taking into account the fact that the elements of
$B_N(\o^{N,1},\ldots, \o^{N,l})$ out of the diagonal are at most
   $N^{-1/2+\eta}=o(1)$ as $N \to \infty$,
       one deduces from (\ref{zlt})
           that for $\delta>0$ small enough
(\ref{base2})  holds true with some $\zeta>0$
   for all $N$ large enough and all $\vec t$ with
   $\|\vec t\|\leq \delta \sqrt{N}$.


\begin{thebibliography}{99}

\bibitem{BFM} H. Bauke, S. Franz, S. Mertens.
   Number partitioning as random energy model. {\it Journal of Stat.
     Mech. : Theory and Experiment,} page P04003 (2004).

\bibitem{BaMe} H. Bauke, S. Mertens. Universality in the level
    statistics of disordered systems. {\it Phys. Rev. E} {\bf 70},
      025102(R) (2004).

\bibitem{BGK} G. Ben Arous, V. Gayrard, A. Kuptsov.
    A new REM conjecture.  Preprint (2006).

\bibitem{BCP} C. Borgs, J. Chayes and B. Pittel.
       Phase transition and finite-size scaling for the integer
partitioning problem. {\it Random Structures and Algorithms} 19,
247-288 (2001).

\bibitem{BCMN1} C. Borgs,  J. Chayes, S. Mertens and C. Nair.
   Proof of the local REM conjecture for number partitioning I:
Constant energy scales. Preprint (2005).
  To appear in {\it Random Structures and Algorithms}

\bibitem{BCMN2} C. Borgs,  J. Chayes, S. Mertens and C. Nair.
 Proof of the local REM conjecture for number partitioning II:
 Growing energy scales. Preprint (2005).
   To appear in {\it Random Structures and Algorithms}

\bibitem{BCMP} C. Borgs, J. Chayes, S. Mertens and B. Pittel.
Phase diagram for the constrained integer partitioning problem. {\it
Random Structures and Algorithms} 24, 315-380 (2004).

\bibitem{BK1} A. Bovier, I. Kurkova. Poisson convergence in the
restricted $k$-partitioning problem. To appear in {\it Random
Structures
  and Algorithms} (2007).

\bibitem{BK2} A. Bovier, I. Kurkova. Local energy statistics in
disordered systems : a proof of the local REM conjecture. {\it
Commun.
  Math. Phys.} {\bf 263}  513--533 (2006).

\bibitem{BK3} A. Bovier, I. Kurkova.
    A tomography of the GREM : beyond the REM conjecture.
{\it Commun.
  Math. Phys.} {\bf 263}  535--552 (2006).

\bibitem{CSY} F. Comets, T. Shiga, N. Yoshida.
  Probabilistic analysis of directed polymers in a random
    environment: a review.
{\it Advanced Studies in Pure Mathematics} {\bf 39} (2004),
 {\it  Stochastic Analysis on Large Scale Interacting Systems,}
   115--142.

\bibitem{D1} B. Derrida. Random-Energy model : an exactly
   solvable model of disordered systems.
     {\it Phys. Rev. B} (3) {\bf 24}(5) 2613--2626(1981).

 \bibitem{D2} B. Derrida. A generalisaton of the random energy
 model that incldes correlations betwen energies. {\it
   Jounal Phys. Lett.}
 {\bf 46}, 401--407(1985).

\bibitem{ET} P. Erdos, S.J. Taylor, Some problems concerning
  the stucture of random walk paths. {\it Acta Math. Acad. Sci.
  Hung.}
    {\bf 11}, 137--162 (1960).

\bibitem{Fe} W. Feller. {\it An Introduction to Probability Theory
  and its Applications.} Volume I.

\bibitem{M1} S. Mertens.  Phase transition
  in the number partitioning problem. {\it Phys. Rev. Lett.}
  {\bf 81}(20), 4281--4284 (1998).

\bibitem{M2} S. Mertens.  A physicist's approach to number
partitioning. {\it Theoret. Comput. Sci.} {\bf 265}(1--2), 79--108,
(2001).


\end{thebibliography}
\end{document}